\documentstyle[11pt,leqno,graphicx,amscd,fancyhdr,amssymb,latexsym,amsbsy,hyperref,xypic,mathrsfs,pstricks,color,verbatim]{amsart}

\include{micro}
\def\da{\delta}
\def\gz{\gamma}
\def\ggz{\Gamma}
\def\az{\alpha}
\def\tz{\theta}
\def\bz{\beta}
\def\ra{\rightarrow}
\def\lra{\longrightarrow}

\def\mc{{\mathcal{C}}}
\def\ms{\mathcal{S}}

\def\dim{\underline{\bf\mbox{dim}}}

\def\t{\tau}

\def\Int{\mbox{Int}}

\def\Hom{\mbox{Hom}}
\def\Ext{\mbox{Ext}}
\def\fs{{\mathfrak{S}}}
\def\ff{{\mathfrak{F}}}
\def\s{{\mathcal{S}}}

\newtheorem{thm}{Theorem}[section]
\newtheorem{lem}[thm]{Lemma}
\newtheorem{cor}[thm]{Corollary}

\newtheorem{prop}[thm]{Proposition}
\newtheorem{rem}[thm]{Remark}

\theoremstyle{definition}
\newtheorem{example}[thm]{Example}

\newtheorem{defn}[thm]{Definition}

\begin{document}

\title{On indecomposable exceptional modules over gentle algebras}
\pagestyle{fancy}
\fancyhead[CO]{On indecomposable exceptional modules over gentle algebras}
\fancyhead[RE]{}
\fancyhead[RO]{\thepage}
\fancyhead[LO]{}
\fancyhead[LE]{\thepage}
\fancyhead[CE]{Jie Zhang}
\lfoot{}
\rfoot{}
\cfoot{}
\renewcommand{\headrulewidth}{0pt}
\renewcommand{\footrulewidth}{0pt}

\author{Jie Zhang}

\address{Department of Mathematical Sciences,
Tsinghua University,
Beijing, P. R. China, 100084} \email{jie.zhang@@math.tsinghua.edu.cn}

\date{\today}


\begin{abstract}We give a characterization of indecomposable exceptional modules over finite dimensional gentle algebras. As an application, we
 study gentle algebras arising from an unpunctured surface and show that a class of indecomposable modules related to
curves without self-intersections, as exceptional modules, are uniquely determined by their dimension vectors.
\end{abstract}
\maketitle


\section{Introduction}
We study in this paper a finite dimensional gentle algebra $A$ over an algebraic closed field $k$. The notion of gentle algebras was first introduced by Assem and Skowro\'nski in \cite{AS87} in order to classify the iterated tilted algebras of type
$\widetilde{A_n}.$  The gentle algebras are special string algebras which were introduced and studied in \cite{BR87, WW85}.
It has been shown in \cite{BR87} that each indecomposable module over a gentle algebra is either a string module $M(w)$ defined by a string $w$, or
a band module $M(b,n,\phi)$ defined by a band $b$, a vector space $V=k^n$ and $\phi\in \mbox{Aut}(V)$.
However, we are only concerned about indecomposable modules without self-extensions in this paper, and band modules which always have self-extensions will be neglected.

We first give a sufficient and necessary condition to determine when two string modules do not have extensions.
To do this, we say two strings $w, v$ have extensions from $w$ to $v$ if $w$ and $v$  can be connected as a string by an arrow from
$w$ to $v$, or $v$ contains a factor string which is also a substring of $w$ or $w^{-1}$,
see precise details in Definition \ref{def-int}. We first present our first theorem which shows that the extensions
of two strings determine the extensions of the corresponding string modules.

\newtheorem{theorem}{Theorem}
\begin{theorem}Let $w,v$ be two strings over a gentle algebra $A$. Then $$\Ext^1_A(M(w),M(v))\neq 0$$ if and only if $w,v$ do not admit extensions
from $w$ to $v.$
\end{theorem}
Note that a necessary condition on $\Ext^1_A(M(w),M(v))\neq 0$ has been given by Schr\"oer in \cite{S99}, see Proposition \ref{yan}.
 We call a module $M$ exceptional if $\Ext^1_A(M,M)=0$. Note that plenty of exceptional modules are constructed in \cite{C11} by using so-called up-and-down graph. A string $w$ is said to have self-extensions if there exist extensions from $w$ to itself.
As a direct application of the above theorem, we give a characterization of
exceptional indecomposable modules over $A$:
\newtheorem*{corollary}{Corollary}
\begin{corollary}Let $M(w)$ be a string module, then $M(w)$ is exceptional if and only if the string $w$ does not admit self-extensions.
\end{corollary}

We give another important application of Theorem 1. We study gentle algebras $G_\ggz$, introduced in \cite{ABCP10,LF09}, arising from triangulations $\ggz$ of an unpunctured surface $(S,M)$. It has been shown in \cite{ABCP10} that there is a bijection $\gz\longleftrightarrow w_\ggz^\gz$ between the curves which are not in $\ggz$ and the strings over $G_\ggz$. We denote by $M^\gz_\ggz=M(w^\gz_\ggz)$ the corresponding string module, see section \ref{def-string}. Under this correspondence, the string modules induced by the curves without self-intersections are indecomposable exceptional modules \cite{BZ11}.
An indecomposable exceptional module $M$ is said to be exceptionally determined by its dimension vector if for any indecomposable exceptional module $N$ $\dim M=\dim N$ implies $M\cong N.$ We show

\begin{theorem}Let $\gz$ be a curve without self-intersections which is not in $\ggz$, then $M_\ggz^\gz$ is exceptionally determined by its dimension vector.
\end{theorem}

Note that to find
indecomposable modules uniquely determined (up to isomorphism) by their dimension vectors is a classical question in representation theory, see \cite{G72} about the proof of the Gabriel theorem, see also recent works \cite{GP12, AD12}. Different from this classical question, we want to
find the indecomposable  modules which are exceptionally determined by their dimension vector. Therefore, $M_\ggz^\gz$ in Theorem 2 might not be determined by its dimension vector which means, there might exist a non-exceptional indecomposable module $N$ such that $\dim M_\ggz^\gz=\dim N,$ see Example \ref{ex-1}. Moreover, Theorem 2 also does not hold for all indecomposable exceptional modules, i.e., there exist two non-isomorphic indecomposable exceptional modules with the same dimension vectors, see Remark \ref{rem}.

However, the  indecomposable modules which are exceptionally determined by their dimension vectors in Theorem 2 are related to many important subjects such as cluster algebras, cluster categories (see \cite{FZ02, BMRRT06, A09, K10}). In fact, these exceptionally determined indecomposable modules together with triangulation $\ggz$
correspond to all the generators in the related cluster algebras, also correspond to
all the objects without self-extensions in
the cluster category $\mc_{(S,M)}$ of the marked surface $(S,M)$ without punctures which has been studied in \cite{A09, BZ11} (we refer to \cite{K10} for more details).

The paper is organized as follows:
we first introduce in section 2 the definition of extensions of two strings
 and present Theorem 1 whose proof will be postponed to section 4.
As an application of Theorem 1, we study, in section 3, the gentle algebras from unpunctured surfaces and give the proof of Theorem 2.
Section 4 is devoted to the proof of Theoreom 1.

\section{Main result}

\subsection{Strings}\label{string-alg}
Recall from \cite{AS87, BR87, S99} that a finite-dimensional algebra $A$ is a
gentle algebra if there is a quiver $Q=(Q_0,Q_1)$ and an admissible ideal $I$
such that  $A=kQ/I$ and the following conditions hold:
\begin{itemize}
\item[$(G1)$]at each vertex of $Q$ start at most two arrows and stop at most two arrows,
\item[$(G2)$]for each arrow $\alpha$ there is at most one arrow $\beta$ and at most one arrow
$\da$ such that $\alpha\beta \not\in I $ and $\da\alpha\not\in I,$
\item[$(G3)$]$I$ is generated by paths of length 2,
\item[$(G4)$]for each arrow $\alpha$ there is at most one arrow $\beta$ and at most one arrow
$\da$ such that $\alpha\beta\in I $ and $\da\alpha\in I.$
\end{itemize}

Note that an algebra $A=kQ/I$ with $I$
an admissible monomial ideal is
called a string algebra if it satisfies the two conditions $(G1)$ and $(G2)$.

For each arrow $\bz$, $s(\bz)$ (resp. $e(\bz)$) denotes its starting point (resp. its ending
point). We denote by $\bz^{-1}$ the formal inverse of $\bz$ with $s(\bz^{-1})=e(\bz)$ and $e(\bz^{-1})=s(\bz)$.
A word $w=\az_1\az_2\cdots\az_n$ of arrows and their formal inverses is called a string of length $n\geq 1$ if $\az_{i+1}\neq \az_i^{-1}, e(\az_i)=s(\az_{i+1})$ for all $1\leq i\leq n-1$, and no subword nor its inverse is in $I$. We denote the length of $w$ by $l(w)$.
 Hence, a string can be viewed as a walk in $Q:$
$$w: e_1 \frac{~~~\az_1~}{}e_2 \frac{~~~\az_2~}{}\cdots e_{n-1} \frac{\az_{n-1}}{}e_{n}\frac{~~~\az_n~}{}e_{n+1}$$
where $e_i\in Q_0$ are vertices of $Q$ and $\alpha_i$ are arrows in either direction. For each vertex $e_i,$ we define
two strings $1_{e^{+1}_i}$, $1_{e^{-1}_i}$ of length 0 with $e(1_{e^t_i})=s(1_{e^t_i})=e_i$ and
$1^-_{e^t_i}=1_{e^{-t}_i}$ for $t\in\{-1,+1\}.$

We denote
by $\mathcal{S}$ the set of all strings over $A$.
A band $b=\az_1\az_2\cdots\az_{n-1}\az_n$ is defined to be a
string $b$ with $s(\az_1)=e(\az_n)$ such that each power $b^m$ is a string, but $b$ itself is not a proper
power of any string.

Recall from \cite{BR87} that each string $w$ defines a unique string
module $M(w)$, each band $b$ yields a family of band modules
$M(b,n,\phi)$ with $n\geq 1$ and $\phi\in \mbox{Aut}(k^n).$  We
refer to \cite{BR87} for more definitions on string modules and band
modules.
\begin{thm}[\cite{BR87}]\label{indec}The indecomposable modules over a string algebra are either string modules or band modules. Moreover,
$\Ext^1_A(M,M)\neq 0$ for each band module $M$.
\end{thm}
\subsection{Maps between string modules}
We present in this subsection a nice basis constructed by Crawley-Boevey for the Hom-space of two string modules \cite{C-B89,S99}.

For a string $w\in\mathcal{S},$ define
${\mathcal{S}}_w=\{(D,E,F)\mid D,E,F\in{\mathcal{S}}, w=DEF\},$
and call $(D,E,F)$ a factor string of $w$ if
 \begin{itemize}
  \item[$(F1)$] $l(D)=0$ or $D=\az_1\az_2\cdots\az_n$ with $\az^{-1}_n\in Q_1,$
  \item[$(F2)$] $l(F)=0$ or $F=\bz_1\bz_2\cdots\bz_m$ with $\bz_1\in Q_1.$
\end{itemize}
Dually, $(D,E,F)$ is said to be a substring of $w$ if the following hold:
 \begin{itemize}
 \item[$(S1)$]  $l(D)=0$ or $D=\az_1\az_2\cdots\az_n$ with $\az_n\in Q_1,$
  \item[$(S1)$]  $l(F)=0$ or $F=\bz_1\az_2\cdots\bz_m$ with $\bz_1^{-1}\in Q_1.$
\end{itemize}
Denote by $\ff(w)$ and $\fs(w)$ the set of all factor strings and substrings of $w$ respectively.
Let $w$ and $v$ be two strings over $A$, a pair $(D, E, F)\times (D', E', F')\in \ff(w)\times\fs(v)$ is said to be an admissible pair (or ad-pair)
if $E=E'$ or $E^{-1}=E'.$ It is easy to understand the ad-pair if one has the following picture in mind.

\medskip

\begin{center}
\includegraphics[height=1.5in]{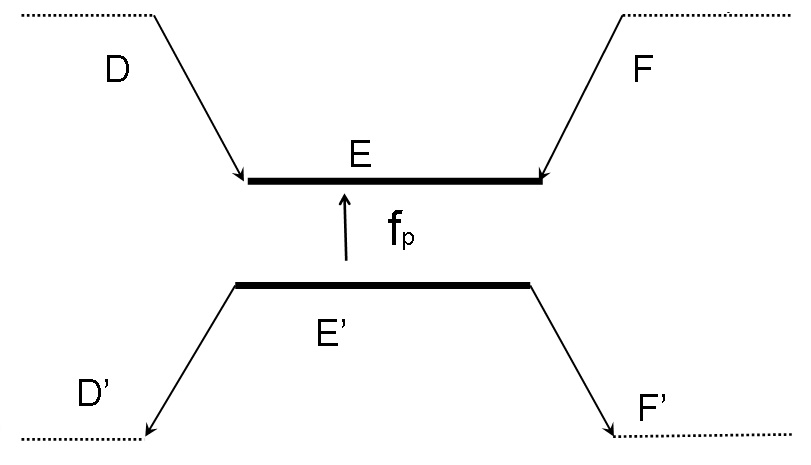}\label{fig}
\end{center}

We call an ad-pair $p=(D,E,F)\times (D',E',F')$ one-sided if one of the following two conditions holds,
\begin{itemize}
  \item[$(O1)$] $l(D)=l(D')=0$ or $l(F)=l(F')=0$ when $E=E'$.
  \item[$(O2)$]  $l(D)=l(F')=0$ or $l(E)=l(F')=0$ when $E^{-1}=E'$.
\end{itemize}
An ad-pair is said to be two-sided if it is not one-sided.

Recall that each ad-pair $p\in \ff(w)\times\fs(v)$ as above yields a canonical homomorphism
$f_p: M(w)\longrightarrow M(v)$ in $\Hom_A(M(w), M(v))$ by identifying the factor module $M(E)$ of $M(w)$ with
the submodule $M(E')$ of $M(v)$. A basis of Hom-space of two string modules can be described as follows:
\begin{thm}[\cite{C-B89}]\label{C-B89} Consider two strings $w$ and $v$ over $A$. Then $\{f_p\mid p\in \ff(w)\times\fs(v)\}$ is a basis of the $\mbox{Hom}_A(M(w),M(v))$.
\end{thm}
\subsection{Extension of two strings}
We present Theorem 1 in this subsection by first introducing the following important definition:
\begin{defn}\label{def-int}
We say two strings $w,v$ have extensions from $w$ to $v$ if one of the following conditions holds
\begin{itemize}
  \item[$(E1)$]There exists $\az\in Q_1$ such that $w\az v$ or $w\az v^{-1}$ is a string.
  \item[$(E2)$]There exists $\bz\in Q_1$ such that $w^{-1}\bz v$ or $w^{-1}\bz v^{-1}$ is a string.
  \item[$(E3)$]$w, v$ admit a two-sided ad-pair in $\ff(v)\times\fs(w)$.
\end{itemize}
\end{defn} A string $w$ is said to have self-extensions if $w=v$ in the above definition.
\begin{example}Consider a gentle algebra $A=kQ/I$ with $Q$ given as follows
\begin{center}
\includegraphics[height=1.5in]{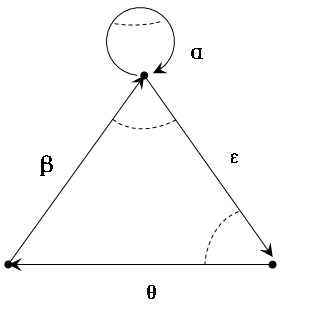}\label{fig}
\end{center}
and $I=<\az^2,\bz\varepsilon,\varepsilon\theta>.$
Then $w=\bz\az^{-1}\bz^{-1}$, $v=\varepsilon^{-1}\az\varepsilon$ have extension from $w$ to $v$ by $(E3).$
In addition, both $\bz$ and $\varepsilon$ have self-extension by $(E1)$ or $(E2)$. However, $\theta$, $\varepsilon$ do not have
extensions from $\varepsilon$ to $\theta,$ note that $\varepsilon$ and $\theta$ only admit an one-sided
ad-pair in $\ff(\theta)\times\fs(\varepsilon).$

\end{example}

Now we state our main theorem whose proof is given in section \ref{pf}
\begin{thm}\label{main-thm}Let $w,v$ be two strings over a gentle algebra $A$. Then $\Ext^1_A(M(w),M(v))\neq 0$ if and only if $w,v$  admit extensions
from $w$ to $v.$
\end{thm}

The following corollary follows directly from the above theorem by taking $v=w.$
\begin{cor}\label{cor-self-ext} Let $w\in\s$, then $M(w)$ is exceptional if and only $w$ does not have self-extensions. That is,
the following conditions hold
\begin{itemize}
  \item[$(1)$]there does not exist $\az\in Q_0$ such that $w\az w\in \s$ or $w\az w^{-1}\in \s,$
  \item[$(2)$]there does not exist $\bz\in Q_0$ such that $w\bz^{-1} w\in \s$ or $w\bz^{-1} w^{-1}\in \s,$
  \item[$(3)$]$w$ only admits one-sided ad-pairs in $\ff(w)\times\fs(w).$
\end{itemize}
\end{cor}

\begin{rem}
  The Theorem \ref{main-thm} does not hold for string algebras as the following example shows. Let $Q$ be the quiver
  $$\circ\overset{\az}{\lra}\circ\overset{\bz}{\lra}\circ\overset{\tz}{\lra}\circ$$
with $I=<\az\bz\tz>.$ Then $A=kQ/I$ is a string algebra and $w=\az\bz, v=\bz\tz$ have extension from $w$ to $v$.
However, $\Ext^1_A(M(w),M(v))=0$ since $M(w)$ is projective.
\end{rem}

\section{An application}
In this section, we study gentle algebras arising from triangulations of a marked surface without punctures. By applying the results in the previous section, we prove the Theorem 2 in the introduction.

\subsection{Triangulation of a marked surface without punctures}\label{notation}
Consider a compact connected oriented 2-dimensional bordered  Riemann surface $S$ and a finite set of marked points $M$ lying on the boundary $\partial S$ of $S$ with at least one marked point on each boundary component.
The condition $M \subset \partial S$ means that we do not allow the marked surface $(S,M)$ to have punctures.

By a curve $\gz$ in $(S,M)$, we mean a continuous function $\gz:[0,1]\rightarrow S$ such that $\gz(0), \gz(1)\in M$ and $\gz(t)\notin M$ for $0<t<1$. An arc is a curve $\gz$ in $(S,M)$ such that $\gz\mid_{(0,1)}$ is injective. We always consider non-contractible curves up to homotopy, and for any collection of curves we implicitly assume that their mutual intersections are minimal possible in their respective homotopy classes. For two curves $\gz',\gz$ in $(S,M)$ we denote by $\Int(\gz',\gz)$ the minimal intersection number of two representatives of the homotopic classes of $\gz'$ and $\gz$ (the intersections at the endpoints do not count).
A triangulation of $(S,M)$ is a maximal
collection $\ggz$ of arcs that do not intersect except at their endpoints.

Let $\ggz=\{\t_1,\t_2,\ldots,\t_n\}$ be a triangulation of $(S,M)$ and $\gz\not\in\ggz$ be a curve  with $d=\sum_{\t\in \ggz}\Int(\t,\gz)$ which is allowed to have self-intersections.
We fix an orientation for the curve $\gz$ and denote $p_0=\gz(0)\in M$ and $p_{d+1}=\gz(1)\in M$. Along its orientation, let $\t_{i_1}, \t_{i_2},\ldots, \t_{i_d}$ be the internal arcs of $\ggz$ that intersect $\gz$ at $p_1=\gz(r_1),\ldots, p_d=\gz(r_d)$ in the fixed orientation of $\gz$, where $0<r_1<r_2<\cdots <r_d<1$. See the following figure.

\begin{center}
\includegraphics[height=1.8in]{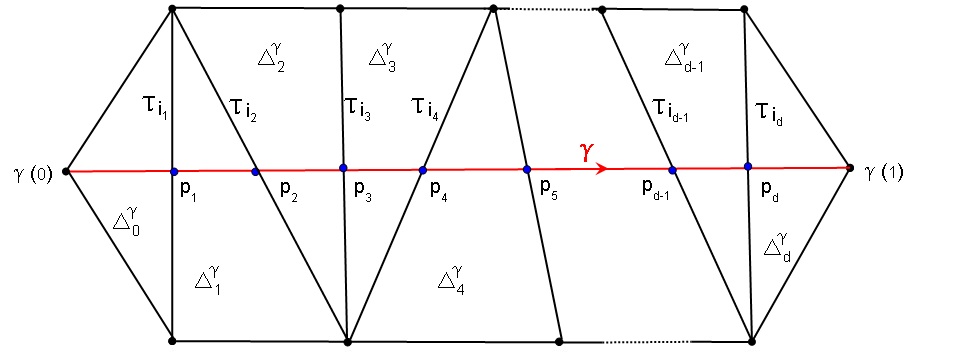}

Figure I.
\end{center}

We denote the restrictions of $\gz$ by $\gz_0={\gz\mid}_{(0,r_1)},$ $\gz'={\gz\mid}_{[r_1,r_d]}$ and $\gz_d={\gz\mid}_{(r_d,1)}.$
Along its way, the curve $\gz$ is passing through (not necessarily distinct) triangles $ \Delta^\gz_0, \Delta^\gz_1, \ldots, \Delta^\gz_d$.

If $\gz$ is an arc which is not in $\ggz$ and $\t^n_j,{\t^n_k}, {\t^n_l}$ are the three arcs of $\Delta^\gz_n$. If  $\Delta^\gz_n\neq\Delta^\gz_d$ or $\Delta^\gz_0$,  then $\Delta^\gz_n\cap\gz$  consists of segments of $\gz$ as follows:

\begin{center}
\includegraphics[height=1.5in]{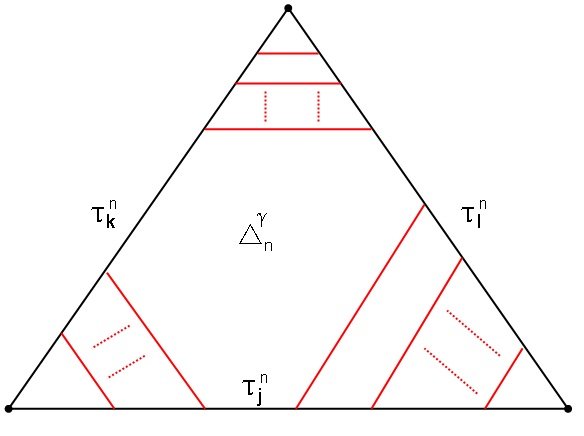}
\end{center}

Note that $\gz$ can cross $\Delta^\gz_n$ many times.
It is easy to see that the sum of intersection numbers of any two arcs with $\gz$ is not less than that of $\gz$ with the third one,
and we have the following:
$$
\left\{
  \begin{array}{ll}
   \Int(\gz,\t^n_k)+\Int(\gz,\t^n_j)\geq \Int(\gz,\t^n_l), & \\
    \Int(\gz,\t^n_k)+\Int(\gz,\t^n_j)+\Int(\gz,\t^n_l)\in 2\mathbb{Z}.\quad\quad\quad\quad\quad\quad & \hbox{(*)}
  \end{array}
\right.
$$

However, if $n=0, d$, then $\Delta^\gz_n\cap\gz$ can be described as follows in two cases:
\begin{center}
\includegraphics[height=1.5in]{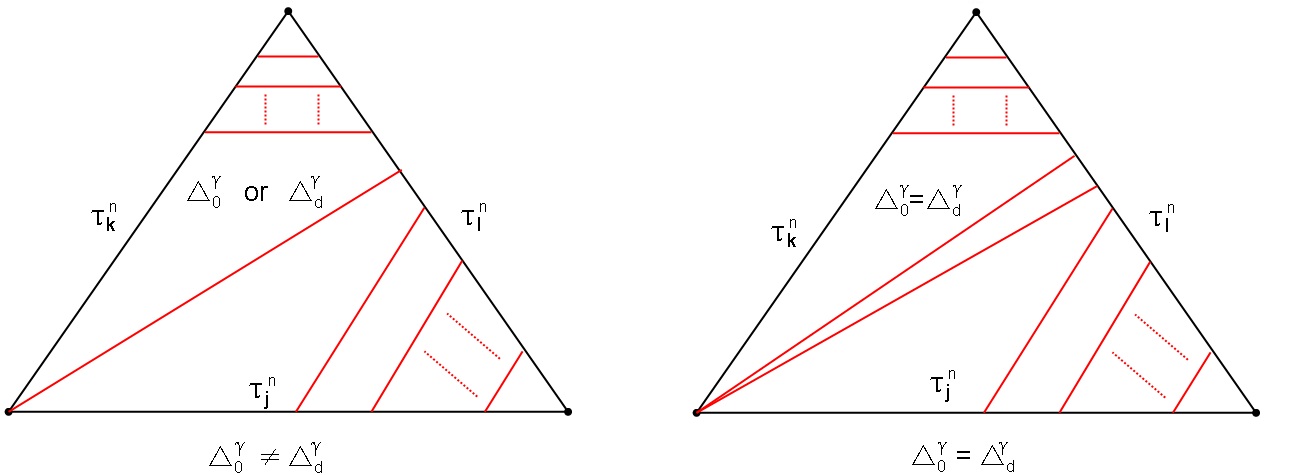}
\end{center}
If $\Delta^\gz_0\neq\Delta^\gz_d,$ then $j,k,l$ can be ordered with $\t^n_l=\t_{i_1}$ or $\t^n_l=\t_{i_d}$ such that
$$\left\{
    \begin{array}{ll}
      \Int(\gz,\t^n_k)+\Int(\gz,\t^n_j)+1=\Int(\gz,\t^n_l), &  \\
     \Int(\gz,\t^n_k)+\Int(\gz,\t^n_j)+\Int(\gz,\t^n_l)\in 2{\mathbb{Z}}+1.\quad\quad\quad\quad\quad & \hbox{(**)}
    \end{array}
  \right.
$$

If $\Delta^\gz_0=\Delta^\gz_d,$ then $j,k,l$ can be ordered with $\t^n_l=\t_{i_1}=\t_{i_d}$ such that
$$\left\{
    \begin{array}{ll}
      \Int(\gz,\t^n_k)+\Int(\gz,\t^n_j)+2=\Int(\gz,\t^n_l), &  \\
     \Int(\gz,\t^n_k)+\Int(\gz,\t^n_j)+\Int(\gz,\t^n_l)\in 2\mathbb{Z}.\quad\quad\quad\quad\quad \quad\quad & \hbox{(***)}
    \end{array}
  \right.
$$
In fact, each arc $\gz$ can also be uniquely recovered from its intersection with each arc in $\ggz$:
\begin{thm}[\cite{M83, T10}]\label{mocher} An arc is uniquely determined by the number of intersections with each
arc of the triangulation.
\end{thm}

\subsection{Gentle algebras arising from $(S,M)$}\label{def-string}
Recall from \cite{DWZ08, LF09, ABCP10} that each triangulation $\ggz$ of $(S,M)$ yields
a quiver with potential $(Q_\ggz,W_\ggz)$:
\begin{itemize}
\item[$(1)$] $Q_{\ggz}=(Q_0,Q_1)$ where the set of vertices $Q_0$ is given by the internal arcs of $\ggz$, and the set of arrows $Q_1$ is defined as follows: Whenever there is a triangle $\vartriangle$ in $\ggz$ containing two
internal arcs $a$ and $b$, then there is an arrow $\rho:a \ra b$ in $Q_1$ if $a$ is a predecessor
of $b$ with respect to clockwise orientation at the joint vertex of $a$ and $b$ in $\vartriangle.$
\bigskip

\item[$(2)$]Every internal triangle $\vartriangle$ in $\ggz$ gives rise to an oriented cycle $\alpha_\vartriangle \beta_\vartriangle\gz_\vartriangle$ in
$Q$, unique up to cyclic permutation of the factors $\alpha_\vartriangle, \beta_\vartriangle, \gz_\vartriangle$. We define
$$W_\ggz =\displaystyle\sum_\vartriangle\alpha_\vartriangle \beta_\vartriangle\gz_\vartriangle $$
where the sum runs over all internal triangles $\vartriangle$ of $\ggz$.
\end{itemize}

From \cite{ABCP10} we know that  the Jacobian algebra $G_\ggz$ of the quiver with potential $(Q_\ggz, W_\ggz)$ is a finite-dimensional gentle algebra, and there is a correspondence
$\gz\longleftrightarrow w^\gz_\ggz$ between curves which are not in $\ggz$ and strings over $G_\ggz:$ For each curve $\gz$ which is not in $\ggz$ with $d=\sum_{\t\in \ggz}I(\t,\gz),$ see Figure I as an example, we obtain a string $w^\gz_\ggz$ in $G_\ggz$:
$$w^\gz_\ggz : \t_{i_1} \frac{~~~\az_1~}{}\t_{i_2} \frac{~~~\az_2~}{}\cdots \t_{i_{d-2}} \frac{\az_{d-2}}{}\t_{i_{d-1}}\frac{~~~\az_{d-1}~}{}\t_{i_d}.$$
We denote by $M_\ggz^\gz$ the corresponding string module $M(w^\gz_\ggz)$ over $G_\ggz.$

We call a curve $\gz$  $\ggz-$rigid if
$\Int(\gz,\gz)\neq 0$ and $\Ext^1_{G_\ggz}(M_\ggz^\gz, M_\ggz^\gz)=0.$ Note that $\ggz-$rigid curves are not arcs by definition.
\begin{example}\label{ex} Consider a triangulation $\ggz$ of an annulus with one marked point on one boundary, two on another:
  \begin{center}
\includegraphics[height=1.5in]{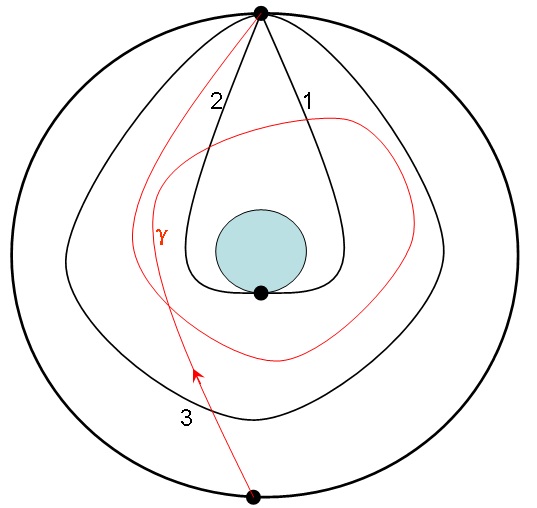}
\end{center}

Then the associated quiver $Q_\ggz$ is as follows:
\begin{center}
\includegraphics[height=1.3in]{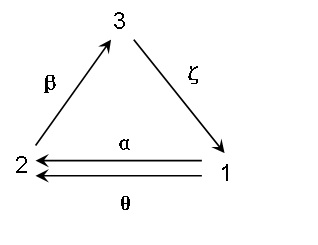}
\end{center}
and $W_\ggz=\az\bz\zeta.$ Let $\gz$ be a curve with self-intersection in the surface as in the above figure, by definition $w^\gz_\ggz=\bz^{-1}\tz$ does not have self-extensions and it is a $\ggz-$rigid curve.
\end{example}

By using Corollary \ref{cor-self-ext}, we give in the following a characterization of indecomposable exceptional modules by curves in $(S,M).$
\begin{lem}If the string $w^\gz_\ggz$ has self-extensions, then $\Int(\gz,\gz)\neq 0.$
\end{lem}
\begin{pf} It suffices to explain Definition \ref{def-int} by the curves in $(S,M).$ Note that there do not exist
loops in $Q_\ggz$.
If there exists $\az\in Q_0$ such that $w^\gz_{\ggz}\az w^\gz_{\ggz}\in \s$, then $\gz$ can be described as follows:
\begin{center}
\includegraphics[height=1.5in]{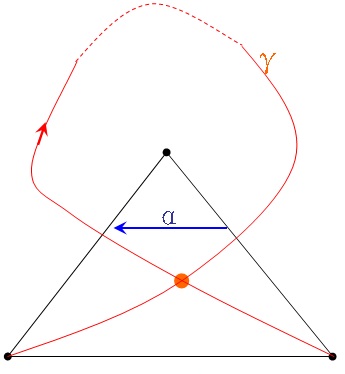}
\end{center}
If there exists $\bz\in Q_0$ such that $ w^\gz_{\ggz}\bz^{-1}  w^\gz_{\ggz}\in \s$,  then $\gz$ can be described similarly as above.

Assume $ w^\gz_{\ggz}$ admits a two-sided ad-pairs
$$(D,E,F)\times(D',E',F')\in\ff( w^\gz_{\ggz})\times\fs( w^\gz_{\ggz}).$$
Then we can characterize $\gz$ as follows when  $E=E'$  or $E^{-1}=E'$ respectively:
\begin{center}
\includegraphics[height=2.5in]{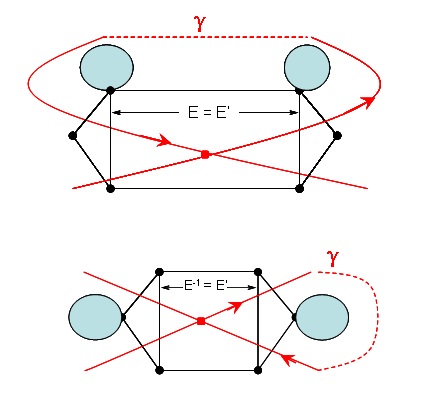}

Figure II
\end{center}
\end{pf}

\begin{cor}[\cite{BZ11}]Let $\gz$ be an arc which is not in $\ggz$, then
$$\Ext^1_{G_\ggz}(M_\ggz^\gz,M_\ggz^\gz)=0$$
for any triangulation $\ggz$ of $(S,M).$
\end{cor}
\begin{pf} Assume that $\ggz$ is a triangulation of $(S,M)$ such that
$\Ext^1_{G_\ggz}(M_\ggz^\gz,M_\ggz^\gz)\neq 0,$ then $w_\ggz^\gz$ has self-extensions by Corollary \ref{cor-self-ext}.
Therefore, $\Int(\gz,\gz)\neq 0$ by the above Lemma which contradicts to the fact that $\gz$ is an arc.
\end{pf}

Combining Theorem \ref{indec} together with the above corollary, indecomposable exceptional modules over $G_\ggz$
are either of the form $M_\ggz^\gz$ with $\gz$ an arc which is not in $\ggz$, or of the form $M_\ggz^\da$ with $\da$ a $\ggz-$rigid curve.

Note that the above corollary does not hold for a $\ggz-$rigid curve $\gz$, that is, there always exists a triangulation $\ggz'$ of $(S,M)$ such that
$\Ext^1_{G_{\ggz'}}(M_{\ggz'}^\gz,M_{\ggz'}^\gz)\neq 0,$ see \cite{BZ11}.
However, we can still find some property for a $\ggz-$rigid curve in terms of Corollary \ref{cor-self-ext}.

Let $\da$ be a $\ggz-$rigid curve, then $\Int(\da,\da)\neq 0$ and $\Ext_{G_\ggz}^1(M^\da_\ggz,M^\da_\ggz)=0.$  We substitute $\da$ for $\gz$ and keep the same notation as in section \ref{notation},
Corollary \ref{cor-self-ext} implies that
\begin{itemize}
  \item[$(\ggz1)$]either $\triangle^\da_0\neq\triangle^\da_n,$ or $\triangle^\da_0=\triangle^\da_n$ with $s(\da)=e(\da).$ This is equivalent to
  the conditions $(1)$ and $(2)$ in Corollary \ref{cor-self-ext},
  \item[$(\ggz2)$]$\da'=\da\mid_{[r_1,r_d]}$ is injective and $\da'$ intersects $\da_0$ or $\da_d$ which means that $w^\ggz_\da$ does not admits two-sided ad-pairs.
\end{itemize}
See the following figure for example:
\begin{center}
\includegraphics[height=2in]{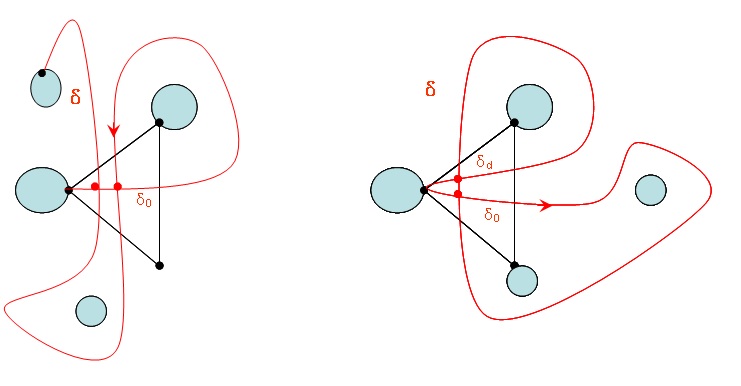}
\end{center}

Roughly speaking, $\da$ only intersects itself with $\da_0$ or with $\da_d.$ However, its restriction $\da'$ is injective. We can do the same discussion
as in section \ref{notation}: Let $\t^n_i, \t^n_j, \t^n_k$ be three arcs of $\Delta^\da_n,$ then  we have similarly as in $(*)$ that
$$\left\{
    \begin{array}{ll}
      \Int(\da,\t^n_k)+\Int(\da,\t^n_j)\geq \Int(\da,\t^n_l),&  \\
      \Int(\da,\t^n_k)+\Int(\da,\t^n_j)+\Int(\da,\t^n_l)\in 2\mathbb{Z}. &
    \end{array}
  \right.
$$
when $\Delta^\da_n\neq\Delta^\da_d$ or $\Delta^\da_0.$
Note that $\da'$ might only intersect one of $\da_0$ and $\da_d$, see the above picture for example. Assume $\da'$ intersects $\da_a$ with $a=0$ or $d$, see the following figure
\begin{center}\label{fig1}
\includegraphics[height=1.5in]{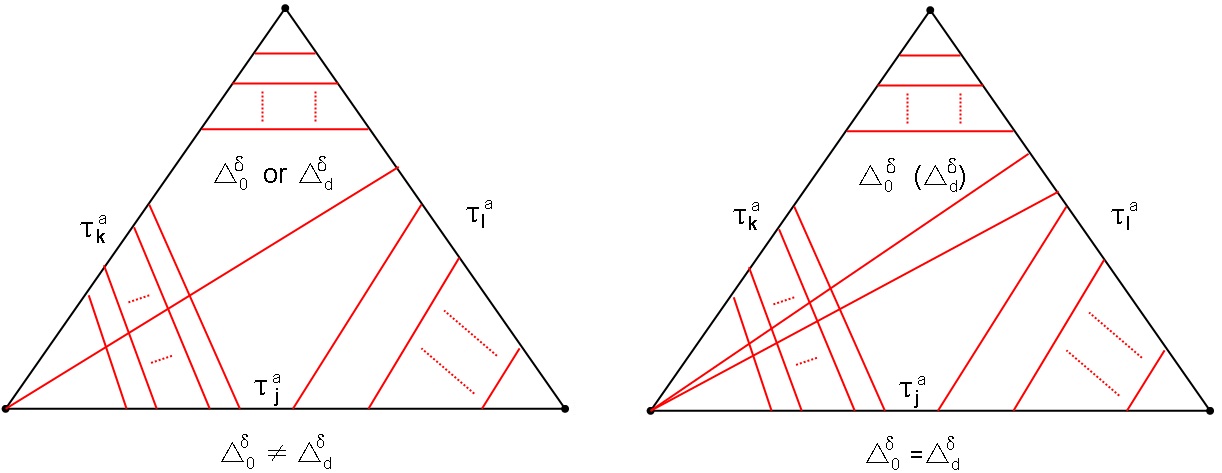}
\end{center}
If $\Delta^\da_0\neq\Delta^\da_d$, then $\Int(\da,\t^a_i)\neq 0$ for $i=j,k,l.$ And
$$\left\{
    \begin{array}{ll}
      \Int(\da,\t^a_k)+\Int(\da,\t^a_j)\geq\Int(\da,\t^a_l)+1, &  \\
      \Int(\da,\t^a_k)+\Int(\da,\t^a_j)+\Int(\da,\t^a_l)\in 2{\mathbb{Z}}+1.\quad\quad\quad\quad & \hbox{(****)}
    \end{array}
  \right.
$$

If $\Delta^\da_0=\Delta^\da_d,$ then $\Int(\da,\t^a_i)\neq 0$ for $i=j,k,l.$ And
$$\left\{\begin{array}{ll}
     \Int(\da,\t^a_k)+\Int(\da,\t^a_j)\geq\Int(\da,\t^a_l), &  \\
     \Int(\da,\t^a_k)+\Int(\da,\t^a_j)+\Int(\da,\t^a_l)\in 2{\mathbb{Z}}. \quad\quad\quad\quad\quad& \hbox{(*****)}
    \end{array}
  \right.$$
By the above observation, we prove in the following the Theorem 2 in the introduction.

\medskip

\emph{{Proof of Theorem 2:}} Let $\ggz=\{\t_1,\ldots,\t_n \}$ be a triangulation of $(S,M)$ and $N$ be an indecomposable exceptional module with
$$\dim M_\ggz^\gz=\dim N.$$
By Theorem \ref{indec}, we can assume that $N=M^\da_\ggz$ with $\da$  either an arc which is not in $\ggz$ or a $\ggz-$rigid curve.

By definition,
$$(\Int(\gz,\t))_{\t\in\ggz}=\dim M_\ggz^\gz=\dim M_{\ggz}^\da=(\Int(\da,\t))_{\t\in\ggz}$$
and $\Int(\da,\t)=\Int(\gz,\t)$
for each $\t\in\ggz.$ Therefore, if $\da$ is an arc which is not in $\ggz$,
then we get the proof by Theorem \ref{mocher}.

We assume now that $\da$ is a $\ggz-$rigid curve with $d=\sum_{\t\in \ggz}\Int(\t,\gz)=\sum_{\t\in \ggz}\Int(\t,\da)$. We are going to prove the theorem by
showing that the dimension vector of $M^\da_\ggz$ (induced by a $\ggz-$rigid curve) never coincides with that of $M^\gz_\ggz$ (induced by an arc which is not in $\ggz$).

Without loss of generality, we suppose that the restriction $\da'$ of $\da$ intersects $\da_0$. Consider $\Delta^\da_0$ with three arcs $\t^0_j,{\t^0_k}, {\t^0_l},$ then $\Int(\da,\t^0_i)\neq 0$ for $i=j,k,l$ by $(****).$ We have the following two cases:
  \begin{itemize}
    \item[(I)]If $\Delta^\da_0\neq \Delta^\da_d$, one has
$$\left\{
  \begin{array}{ll}
   \Int(\da,\t^0_k)+\Int(\da,\t^0_j)\geq\Int(\da,\t^0_l)+1, & \\
   \Int(\da,\t^0_k)+\Int(\da,\t^0_j)+\Int(\da,\t^0_l)\in 2{\mathbb{Z}}+1, &
  \end{array}
\right.$$
Therefore, there must exist a triangle $\Delta^\gz_m=\Delta^\da_0$  with the same three arcs $\t^0_j,{\t^0_k}, {\t^0_l}$ satisfying $\Int (\gz, \t^0_i)\neq 0$ for $i=j,k,l,$ and
\begin{itemize}
  \item[$\bullet$]$\Int(\gz,\t^0_k)+\Int(\gz,\t^0_j)\geq\Int(\gz,\t^0_l)+1$ which contradicts to $(**),$ since $(**)$ implies that there exists
$\{\t'_j,\t'_k,\t'_l\}=\{\t^0_j,\t^0_k,\t^0_l\}$ such that
$$\Int(\gz,\t'_k)+\Int(\gz,\t'_j)<\Int(\gz,\t'_l).$$
  \item[$\bullet$]$\Int(\gz,\t^0_k)+\Int(\gz,\t^0_j)+\Int(\gz,\t^0_l)\in 2{\mathbb{Z}}+1$ which contradicts to $(*),(***).$
\end{itemize}
\item[(II)]If $\Delta^\da_0=\Delta^\da_d$, then we have by $(*****)$ that, for each $0\leq n\leq d$
$$\left\{
    \begin{array}{ll}
     \Int(\da,\t^n_k)+\Int(\da,\t^n_j)\geq\Int(\da,\t^n_l), & \\
      \Int(\da,\t^n_k)+\Int(\da,\t^n_j)+\Int(\da,\t^n_l)\in 2{\mathbb{Z}}, &
    \end{array}
  \right.$$
where $\t^n_j, \t^n_k, \t^n_l$ are three arcs of triangle $\Delta^\da_n.$ Therefore, we also have in $\Delta^\gz_0$ with three arcs $\t'_k,\t'_j,\t'_l$ that
$$\left\{
    \begin{array}{ll}
     \Int(\gz,\t'_k)+\Int(\gz,\t'_j)\geq\Int(\gz,\t'_l), & \\
     \Int(\gz,\t'_k)+\Int(\gz,\t'_j)+\Int(\gz,\t'_l)\in 2{\mathbb{Z}}.&
    \end{array}
  \right.$$
We consider the following two subcases:
   \begin{itemize}
     \item[(II.1)]If $\Delta^\gz_0\neq\Delta^\gz_d,$ it contradicts to $(**).$

     \item[(II.2)]If $\Delta^\gz_0=\Delta^\gz_d,$ it contradicts to $(***)$ which implies that $j,k,l$ can be ordered such that
$$\Int(\gz,\t'_k)+\Int(\gz,\t'_j)<\Int(\gz,\t'_l).$$
   \end{itemize}
\end{itemize}
Therefore, the dimension vector of $M^\da_\ggz$ is always different from that of $M^\gz_\ggz$ which
completes the proof. \quad\quad\quad\quad\quad\quad\quad\quad\quad\quad\quad\quad\quad\quad\quad\quad\quad\quad\quad\quad $\Box$

\medskip

Note that in the proof of the above theorem, we do need that $N$ is an exceptional module. If not, the Theorem 2 does not hold, see the following example:
\begin{example}\label{ex-1}Consider the following annulus with a triangulation $\ggz$:
\begin{center}
\includegraphics[height=2in]{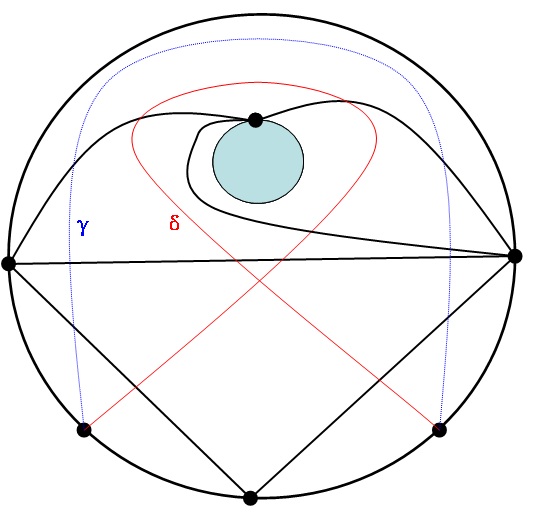}
\end{center}
where $\gz$ is an arc and $\da$ is a curve with self-intersection. It is easy to find that $M^\da_\ggz$ has self-extensions by Theorem \ref{main-thm}. However,
$\dim M^\da_\ggz=\dim M^\gz_\ggz.$
\end{example}
\begin{rem}\label{rem}The above theorem does not hold in general for any indecomposable exceptional module.
Reconsider Example \ref{ex}, let $\da$ be a $\ggz-$rigid curve as follows
\begin{center}
\includegraphics[height=1.5in]{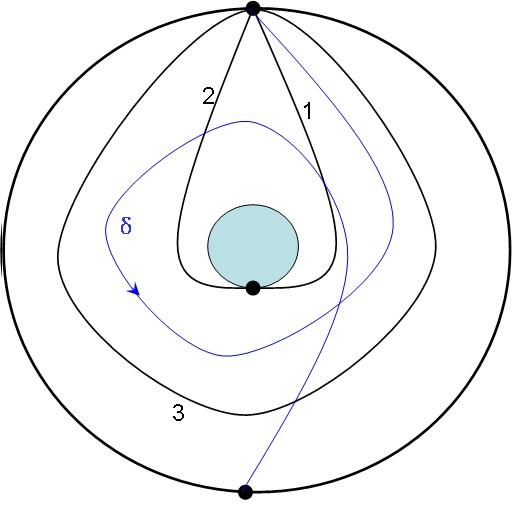}
\end{center}
it is easy to find that $\dim M_\ggz^\gz=\dim M_{\ggz}^\da$ and $\Ext^1_{G_\ggz}(M_\ggz^\gz,M_\ggz^\gz)=0=\Ext^1_{G_\ggz}(M_{\ggz}^\da,M_{\ggz}^\da).$ However, $\gz$ is not homotopic to $\da$ which implies that $M^\gz_\ggz\ncong N^\da_\ggz.$
\end{rem}


\section{Proof of the Theorem \ref{main-thm}}\label{pf}

We prove in this section the Theorem \ref{main-thm} by Auslander-Reiten formula in a combinatorial way. Let first recall from \cite{BR87} the Auslander-Reiten theory for a gentle algebra.
\subsection{Auslander-Reiten theory over gentle algebras}
 Let $A=kQ/I$ be a finite-dimensional gentle algebra with $Q=(Q_0,Q_1)$ and $\mathcal{S}$ the set of all strings over $A$.
 Similarly as in \cite{BR87,S99}, one can define
two maps from $Q_1$ to $\{1,-1\}$ to deal with the concatenation of two strings, especially for the strings of length zero.
However, in this paper the strings of length zero can always be chosen such that the concatenation is well-defined and we only use $1_{e_i}$ to denote
the chosen one, i.e., either $1_{e^{-1}_i}$ or $1_{e^{+1}_i}.$

As for a string $w=\az_1\az_2\cdots\az_n$, we denote by $e(w)=e(\az_n),$ and $s(w)=s(\az_1).$ We say $w$ starts (or ends)  on a peak if there is no arrow $\az \in Q_1$ with $\az w \in \ms$ (or $ w\az^{-1} \in \ms$); likewise, a string $w$ starts (or ends) in a deap if there is no arrow $\bz \in Q_1$ with $\bz^{-1}w \in \ms$ (or $ w\bz \in \ms$).

A string $w=\az_1\az_2\cdots\az_n$ with all $\az_i \in Q_1$ is called direct string, and a string of the form $w^{-1}$ where $w$ is a direct string is called inverse string.
Strings of length zero are both direct and inverse. For each arrow $\az \in Q_1$, let $N_\az= {V_\az}\az U_\az$ be the unique string such that $U_\az$ and $V_\az$ are inverse strings and $N_\az$ starts in a deep and ends on a peak (see the following figure).
\begin{center}
\includegraphics[height=.25in]{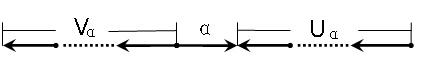}
\end{center}
If the string $w$ does not start on a peak, we define $_hw=V_\az\az w$ and say that $_hw$ is obtained from $w$ by adding a hook on the starting point $s(w)$. Dually, if $w$ does not end on a peak , we define $w_h=w\bz^{-1}{V_\bz}^{-1}$ and say that $w_h$ is obtained from $w$ by adding a hook on the ending point $e(w)$ (see the following figure).
\begin{center}
\includegraphics[height=1in]{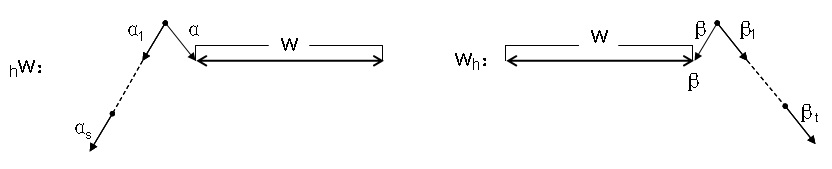}
\end{center}
Suppose now that the string $w$ starts on a peak. If $w$ is not a direct string, we can write $w=\az_1\cdots\az_s\az^{-1} {_cw}=U_\az^{-1}\az^{-1} {_cw}$ for some $\az\in Q_1$ and $s\geq 0$. We say in this case that $_cw$ is obtained from $w$ by deleting a cohook on $s(w)$. If $w$ is a direct string, we define $_cw=0$.
Dually, assume that $w$ ends on a peak. Then, if $w$ is not an inverse string, we can write $w=w_c\bz\bz_1^{-1}\cdots\bz_t^{-1}=w_c\bz U_\bz$ for some $\bz\in Q_1$ and $t\geq 0$. We say that $w_c$ is obtained from $w$ by deleting a cohook on $e(w)$, and if $w$ is an inverse string, we define $w_c=0$ (see the following figure).

\begin{center}
\includegraphics[height=1.2in]{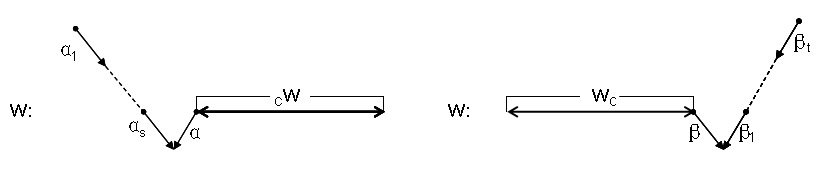}
\end{center}

The following theorem summarizes the results from  \cite{BR87} concerning AR-sequences between string modules:

\begin{thm}[\cite{BR87}]\label{BR87} For  a string algebra $A$, let $w$ be a string such that  $M(w)$ is not an injective $A-$module. Then the AR-sequence starting in $M(w)$ is given
\begin{itemize}
\item[(1)] if $w$ neither starts nor ends on a peak by
 $$0\lra M(w) \lra M(w_h)\oplus M(_hw) \lra M(_hw_h) \lra 0,$$

\item[(2)] if $w$ does not start but ends on a peak by
$$0\lra M(w) \lra M(w_h)\oplus M(_cw) \lra M(_cw_h) \lra 0,$$

\item[(3)] if $w$ starts but does not end on a peak by
$$0\lra M(w) \lra M(w_c)\oplus M(_hw) \lra M(_hw_c) \lra 0,$$

\item[(4)] if $w$ both starts and ends on a peak by
$$0\lra M(w) \lra M(w_c)\oplus M(_cw) \lra M(_cw_c) \lra 0.$$
\end{itemize}
\end{thm}
\subsection{Proof of Theorem \ref{main-thm}}
We first consider the special case: $M(v)$ is an injective module or $M(w)$ is an projective module.
\begin{lem}\label{lem-inj} Let $w,v$ be two strings over $A$,
if either $M(v)$ is injective or $M(w)$ is projective, then $w,v$ do not have extensions from $w$ to $v$
\end{lem}
\begin{pf}We only consider the case that $M(v)$ is injective which means there does not exist an arrow $\az$ or $\bz$ such that $\az v$ or $v\bz^{-1}$ is a string. Hence
$(E1)$ and $(E2)$ do not hold for any $w$ with $v.$

Assume $p=(D,E,F)\times (D',E',F')\in\ff(v)\times\fs(w)$, then either $l(D)=0$ or $l(F)=0$ by property of injective modules. Without loss of generality,
we let $l(D)=0$ and consider the following two cases:
\begin{itemize}
\item[$(1)$] If $E=E'$, then $l(D')=0.$ Otherwise, there exists an arrow $\az$ such that $\az E'=\az E$ is a string by definition of substrings. It contradicts to
the fact that $v$ starts on a peak.
\item[$(2)$] If $E^{-1}=E'$, then $l(F')=0$. Otherwise there exists an arrow $\bz$ such that $E'\bz^{-1}=(\bz E)^{-1}$ is a string by definition of substrings. It contradicts again to
the fact that $v$ starts on a peak.
\end{itemize}
Therefore any $p\in\ff(v)\times\fs(w)$ is an one-sided ad-pair. This completes the proof.
\end{pf}

When $w,v$ admit a two-sided ad-pair in $\ff(v)\times\fs(w),$ Schr\"oer constructed in Proposition 4.9 of \cite{S99} a non-split exact sequence
from $M(v)$ to $M(w)$ which gives the following proposition:

\begin{prop}[\cite{S99}]\label{yan} If $w, v$ admit a two-sided ad-pair in $\ff(v)\times\fs(w)$, then
$$\Ext^1_A(M(w),M(v))\neq 0.$$
\end{prop}
The following corollary comes directly from Definition \ref{def-int} and the above proposition.
\begin{cor}\label{cor-imp}Let $w, v$ be two strings which have extensions from $w$ to $v,$ then
$$\Ext^1_A(M(w),M(v))\neq 0.$$
\end{cor}
\begin{pf} In case $(E1)$, let $u=w\az v$ (or $u=w\az v^{-1}$) if $w\az v$ (or $w\az v^{-1}$) is a string.
Then there is a non-split exact sequence
$$0\lra M(v)\lra M(u)\lra M(w)\lra 0$$
which implies that $\Ext^1_A(M(w),M(v))\neq 0.$
The case $(E2)$ is similar to the case $(E1)$, and the case $(E3)$ follows from Proposition \ref{yan}.
\end{pf}

To prove the main theorem, we need to find out that the morphisms induced by some ad-pairs factor through projective modules.

\begin{lem}\label{lem-important}Let $w, v$ be two strings over $A$ and suppose $v=V_\az \az v_0$ for some $\az\in Q_1, v_0\in {\mathcal{S}}.$ If $p=(D,E,F)\times (D',E',F')\in\ff(v)\times\fs(w)$ such that
$DE=V_\az$ with $l(E)>0$. Then $f_p$ factors through a projective module if one of the following conditions holds
\begin{itemize}
  \item[(1)]$E=E'$ and $l(F')>0.$
  \item[(2)]$E^{-1}=E'$ and $l(D')>0.$
\end{itemize}
\end{lem}
\begin{pf} We only prove $(1)$ since $(2)$ follows from $(1)$ if we consider $w^{-1}$ instead of $w.$

Let $V_\az=\az_n^{-1}\cdots\az^{-1}_2\az^{-1}_1$ with $n\geq 1$ since $l(V_\az)\geq l(E)>0.$
Suppose $D=\az^{-1}_n\az^{-1}_{n-1}\cdots\az^{-1}_{j+1}$, $E=\az^{-1}_{j}\cdots\az^{-1}_2\az^{-1}_1,$ where $1\leq j\leq n$ with
$\az^{-1}_{n+1}:=1_{e(\az_n)}.$

By definition of substring, $l(F')>0$ implies that there exists
an arrow $\bz$ such that $w=D'E'F'=D'E'\bz^{-1}F'_0$ with $e(\bz)=s(\az), F'_0\in {\mathcal{S}}$  Since $E'\bz^{-1}=E\bz^{-1}$ is a string, $\bz\az_1\not\in I$. It means
$\bz\az\in I$ by definition of gentle algebras. Let $u=V_\az\bz^{-1}V^{-1}_\bz=DE\bz^{-1}\bz_1\bz_2\cdots\bz_m$ with $m\geq 0$ and $\bz_0:=1_{s(\bz)}$, then $P_{s(\bz)}=M(u)$ is the indecomposable projective module corresponding to $s(\bz)\in Q_0$. One gets the following one-sided ad-pair:
  $$b=(1_{e(\az_n)}, V_\az, F)\times (1_{e(\az_n)}, V_\az, \bz^{-1}V^{-1}_\bz)\quad \quad \quad  (\vartriangle)$$
  $$\in \ff(v)\times \fs(u).$$
  See the following figure
\begin{center}
\includegraphics[height=1.5in]{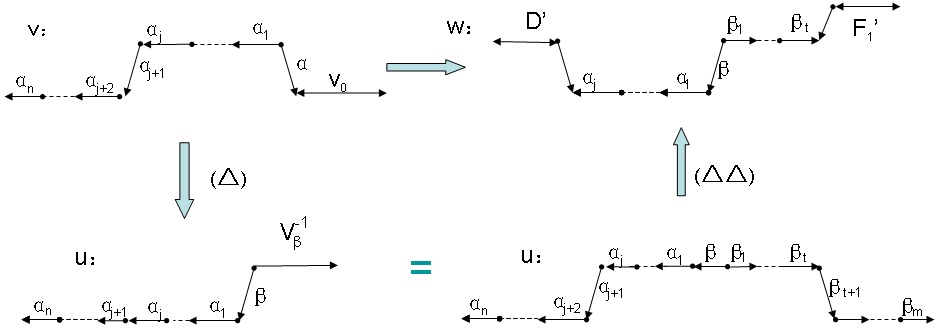}\label{fig}
\end{center}
  On the other hand, by definition of $U_\bz$ one can rewrite
  $w=D'E'F'=D'E'\bz^{-1}\bz_1\bz_2\cdots\bz_t F'_1$ where $0\leq t\leq m, F'_1\in\s$ such that
  $$(D', E'\bz^{-1}\bz_1\cdots\bz_t, F'_1)\in \fs(w).$$
Therefore, we get an ad-pair
  $$c=(\az^{-1}_n\cdots\az^{-1}_{j+1}, E'\bz^{-1}\bz_1\cdots\bz_t, \bz_{t+1}\cdots\bz_m)\times(D', E'\bz^{-1}\bz_1\cdots\bz_t, F'_1)\quad \quad (\vartriangle\vartriangle)$$
  $$\in\ff(u)\times\fs(w)$$
See the above figure.
Combine $(\vartriangle)$ together with $(\vartriangle\vartriangle)$, we get $f_p=f_bf_c$ which completes the proof.
\end{pf}

The following corollary follows directly from the above lemma if we consider $v^{-1}$ in stead of $v.$
\begin{cor}Let $w, v$ be two strings and suppose $v=v_0\bz U_\bz$ for some $\bz\in Q_1, v_0\in {\mathcal{S}}.$ If $b=(D,E,F)\times (D',E',F')\in\ff(v)\times\fs(w)$ such that $EF=V^{-1}_\bz$ with $l(E)>0$. Then $f_b$ factors through a projective module if one of the following conditions holds
\begin{itemize}
  \item[(1)]$E=E'$ and $l(D')>0.$
  \item[(2)]$E^{-1}=E'$ and $l(F')>0.$
\end{itemize}
\end{cor}

In the following, we always assume that $w,v$ do not have extensions from $w$ to $v$, and show that $\Ext^1_A(M(w),M(v))=0$ by considering $v$ in four different cases as in Theorem \ref{BR87}.

\begin{lem}\label{case1}If $v$ does not start but ends on a peak, then
$$\Ext^1_A(M(w),M(v))=0.$$
\end{lem}
\begin{pf}By Auslander-Reiten formula and Theorem \ref{BR87}, we have
$$\Ext_A^1(M(w),M(v))\cong\underline{\Hom}_A(M(_hv_c),M(w)).$$
By Theorem \ref{C-B89}, it suffices to show that that each $f_p\in \Hom_A(M(_hv_c),M(w))$ with $p=(D,E,F)\times (D',E',F')\in \ff({_hv}_c)\times \fs(w)$
factors through a projective module.

Since $v$ does not start but ends on a peak, there exists $\az\in Q_1$ such that $\az v$ is a string.

\begin{itemize}
  \item[(1)]If $v=U_\az$ is an inverse string, then
$$\t^{-1}M(v)=M({(_hv)}_c)=M(V_\az).$$
   \item[(2)]If $v$ is not an inverse string, one can rewrite
$v=v_0\bz\bz^{-1}_1\cdots\bz^{-1}_m=v_0\bz U_\bz$ with $\bz\in Q_1, v_0\in {\mathfrak{S}}$, where $m\geq 0$, $\bz_0:=1_{e(\bz)}.$
By Theorem \ref{BR87}
$$\t^{-1}M(v)=M({_hv}_c)=M(V_\az \az v_0).$$
  \end{itemize}

{\emph{Claim:}} One gets for each $p=(D,E,F)\times (D',E',F')\in \ff({_hv}_c)\times \fs(w)$ that
$$DE=V_\az.$$

Assume $DE\neq V_\az,$ we are going to prove the claim by constructing a two-sided ad-pair in $\ff(v)\times\fs(w),$ hence get a contradiction
to the assumption.

If $v=U_\az$ is an inverse string, then ${_hv}_c=V_\az$, by definition of  factor string, there is nothing to prove. i,e., one has
$DE=V_\az, F=1_{s(\az)}.$

Assume now that $v$ is not an inverse string, we write $v=v_0\bz U_\bz$, then ${_hv}_c=V_\az\az v_0$. By definition of factor string, $l(DE)\geq l(V_\az)$ and $D\neq V_\az\az.$ We suppose that $l(DE)>l(V_\az)$ and consider the following two cases:
\begin{itemize}
  \item[I.] $l(D)>l(V_\az\az).$ Assume $D=V_\az\az D_0$ with $D_0\in {\mathcal{S}}$ such that $l(D_0)>0$, then
 $$(D_0, E, F\bz U_\bz)\in \ff(v).$$
 Hence there exists a two-sided ad-pair
   $$(D_0, E, F\bz U_\bz)\times (D',E',F')\in \ff(v)\times \fs(w),$$
since $l(D_0)>0$ and $l(\bz U_\bz)>0.$
See the following figure.
   \begin{center}
\includegraphics[height=2in]{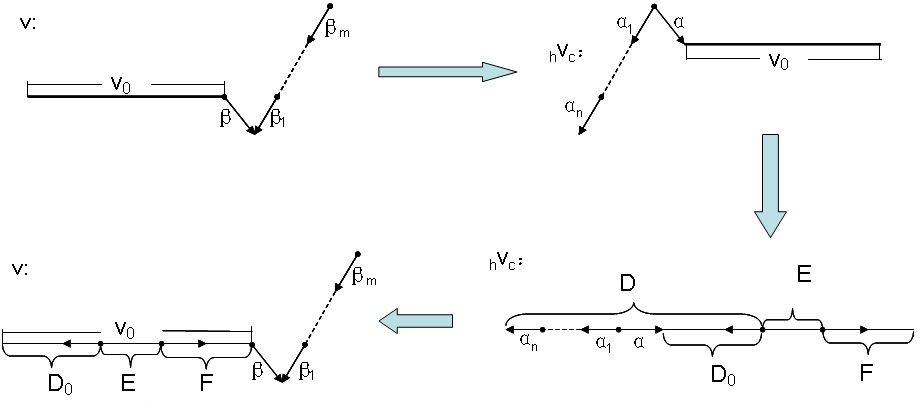}\label{fig}
\end{center}
It contradicts to our assumption.
 \item[II.] $l(D)< l(V_\az\az).$ Let $V_\az=\az^{-1}_n\cdots\az^{-1}_2\az^{-1}_1$ with $n\geq 0$ and $\az^{-1}_0:=1_{s(\az)},$
 $D=\az^{-1}_n\az^{-1}_{n-1}\cdots\az^{-1}_{j+1}$ where $0\leq j\leq n$ and $\az^{-1}_{n+1}:=1_{e(\az_{n})}.$
 We can write $E=\az^{-1}_{j}\cdots\az^{-1}_2\az^{-1}_1\az E_0$  with $E_0\in{\mathcal{S}}$ such that $l(E_0)\geq0$ since $l(DE)>l(V_\az)$.
  Then $$(1_{e(\az)}, E_0, F\bz U_\bz)\in \ff(v).$$
See the following figure.
     \begin{center}
\includegraphics[height=2in]{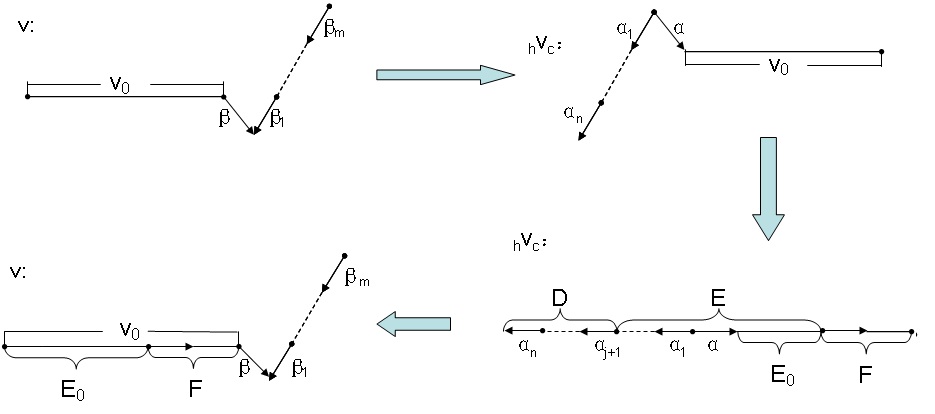}\label{fig}
\end{center}
If $E'=E=\az^{-1}_{j}\cdots\az^{-1}_2\az^{-1}_1\az E_0$, then
  $$(1_{e(\az)}, E_0, F\bz U_\bz)\times (D'\az^{-1}_{j}\cdots\az^{-1}_2\az^{-1}_1\az, E_0, F')\in \ff(v)\times\fs(w)$$
is a two-sided ad-pair since $l(D'\az^{-1}_{j}\cdots\az^{-1}_1\az)>0, l(F\bz U_\bz)>0.$

Similarly, if $E'=E^{-1}$, one gets the following two-sided ad-pair
   $$(1_{e(\az)}, E_0, F\bz U_\bz)\times (D', E^{-1}_0,\az^{-1}\az_1\cdots\az_{j}F')\in \ff(v)\times\fs(w).$$
In both cases, we get a contradiction to the assumption.
  \end{itemize}
It completes the proof of the claim i.e., $$DE=V_\az.$$
Then we use Lemma \ref{lem-important} to show that $f_p$ factors through a projective module.
We first consider the following two cases when $l(E)=l(E')>0:$
\begin{itemize}
  \item[$(1)$]Assume $E=E'$, then $l(F')>0$. Otherwise $e(w)=e(E')=e(E)=s(\az)$, $E\az$ is a string implies
  that $w\az$ is also a string. The definition of gentle algebras  guarantees that $w\az v$ is a string which contradicts
  to the assumption that $w, v$ do not have extensions from $w$ to $v$. Therefore, $l(F')>0.$
   By Lemma \ref{lem-important}, $f_p$ factors through
  projective modules.
  \item[$(2)$]If $E^{-1}=E'$, then $l(D')>0$. Otherwise $s(w)=s(E')=e(E)=s(\az)$, $E\az$ is a string implies ${E'}^{-1}\az$ is a string. Hence
$w^{-1}\az v$ is a string which contradicts
  to the assumption. This implies that $l(D')>0$. Similarly, Lemma \ref{lem-important} implies that $f_p$ factors through
  projective modules.
\end{itemize}
Now, assume $l(E)=l(E')=0$, then $l(F')>0$ or $l(D')>0$. Otherwise, $w=D'E'F'=1_{s(\az)^t}$ with $t=1,-1.$ Hence, either $w\az v$ or $w^{-1}\az v$ is a string which is a contradiction to the hypothesis.

Without loss of generality, let $l(F')>0$ which means there exists
$\bz\in Q_1$ with $e(\bz)=s(\az)$. We rewrite $w=D'E'F'=D'E'\bz^{-1}F'_0$ with $F'_0\in {\mathcal{S}}$ and consider the following two cases:

\begin{itemize}
  \item[(1)]If $\bz\az\in I$, we consider the indecomposable projective module $P_{s(\bz)}$ corresponding to $s(\bz)\in Q_0$, then one gets $f_p$ factors through $P_{s(\bz)}$
as in Lemma \ref{lem-important}.
\item[(2)]If $\bz\az\not\in I$, then $l(D')>0$. Otherwise $w=\bz^{-1}F'_0.$ By definition of gentle algebra, we get $v^{-1}\az^{-1}\bz^{-1}F'_0=v^{-1}\az^{-1}w$ is a string which is again a contradiction.  Therefore, $l(D')>0$ and there exists $\da\in Q_1$ such that $e(\da)=s(\az)$ with $\da\az\in I.$ Similarly as in Lemma \ref{lem-important}, one can find that $f_p$ factors through the indecomposable projective module $P_{s(\da)}$ corresponding to $s(\da)\in Q_0$.
\end{itemize}
Therefore by Auslander-Reiten formula, one gets
$$\Ext^1_A(M(w),M(v))=0.$$
\end{pf}

The following corollary comes directly from the above lemma by considering $v^{-1}$ instead of $v.$
\begin{cor}\label{case2}If $v$ does not end but starts on a peak, then
$$\Ext^1_A(M(w),M(v))=0.$$
\end{cor}

\begin{lem}\label{case3}If $v$ neither starts nor ends on a peak, then
$$\Ext^1_A(M(w),M(v))=0.$$
\end{lem}
\begin{pf}Assume $v$ neither starts nor ends on a peak, then there exists $\az,\bz\in Q_1$ such that $\az v\bz^{-1}$ is a string. Let $u={_hv}_h=V_\az\az v\bz^{-1}V^{-1}_\bz$,
Theorem \ref{BR87} yields that $$\t^{-1}M(v)=M(u).$$
We prove that each non-zero morphism $f_p\in\Hom_A(\t^{-1}M(v),M(w))$ factors through projective modules, where
$p=(D,E,F)\times (D',E',F')\in \ff(u)\times \fs(w)$.

\emph{Claim:} Either $DE=V_\az \mbox{~~or~~} EF=V^{-1}_\bz.$

Similarly as the proof of Lemma \ref{case1}, assume $DE\neq V_\az, EF\neq V^{-1}_\bz$ and we are going to prove the claim by constructing a two-sided ad-pair
in $\ff(v)\times\fs(w)$ which contradicts to our assumption.

By definition of factor string, we have
\begin{itemize}
  \item[(1)]$l(DE)\geq l(V_\az)$ with $l(D)\leq l({_hv} \bz^{-1})$ and $D\neq V_\az\az,$
  \item[(2)]$l(EF)\geq l(V^{-1}_\bz)$ with $l(F)\leq l(\az {v_h})$ and $F\neq \bz^{-1}V^{-1}_\bz.$
\end{itemize}
If $DE=V_\az$ or $D={_hv} \bz^{-1}$, then there is nothing to prove for the claim.
We assume that $l(DE)>l(V_\az)$ with $D\neq {_hv} \bz^{-1}$ and consider the following two cases:

\begin{itemize}
  \item[I.] $l(V_\az\az)<l(D)<l({_hv} \bz^{-1}).$ We rewrite $D=V_\az\az D_0$ with $D_0\in{\mathcal{S}}$ such that $l(D_0)>0$. Note that
  $l(F)<l(\az {v_h})$ and $F\neq \bz^{-1}V^{-1}_\bz$ by $(2)$ above,
   we consider following two subcases related to $l(F):$
  \begin{itemize}
    \item[(I.1)]$l(\az {v_h})>l(F)>l(\bz^{-1}V^{-1}_\bz).$ We rewrite $F=F_0\bz^{-1}V^{-1}_\bz$ with $F_0\in{\mathcal{S}}$ such that $l(F_0)>0.$ Then there exists an ad-pair
   $$(D_0, E, F_0)\times (D',E',F')\in \ff(v)\times \fs(w)$$ which is not one-sided since $l(D_0)>0$ and $l(F_0)>0.$ See the following figure
        \begin{center}
\includegraphics[height=2in]{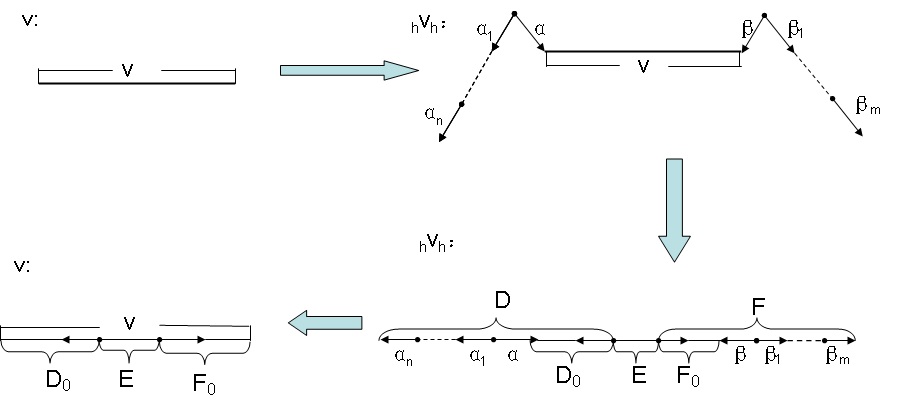}\label{fig}
\end{center}

    \item[(I.2)]$l(F)\leq l(V^{-1}_\bz).$ Let $V^{-1}_\bz=\bz_1\bz_2\cdots\bz_m$ with $m\geq 0$ where $\bz_0:=1_{s(\bz)}$, and rewrite         $F=\bz_{t+1}\cdots\bz_{m-1}\bz_m$ with $0\leq t\leq m$ where $\bz_{m+1}:=1_{e(\bz_m)}.$ We can rewrite $E=E_0\bz^{-1}\bz_1\cdots\bz_{t}$ with $E_0\in{\mathcal{S}}$. Then we obtain a two-sided ad-pair
    $$(D_0, E_0, 1_{e(\bz)})\times (D', E_0, \bz^{-1}\bz_1\cdots\bz_{t}F')$$
    $$\in\ff(v)\times \fs(w)$$
  when $E=E'$, since $l(D_0)>0$ and $l(\bz^{-1}\bz_1\cdots\bz_{t-1}F')>0$.
  See the following figure
\begin{center}
\includegraphics[height=1.7in]{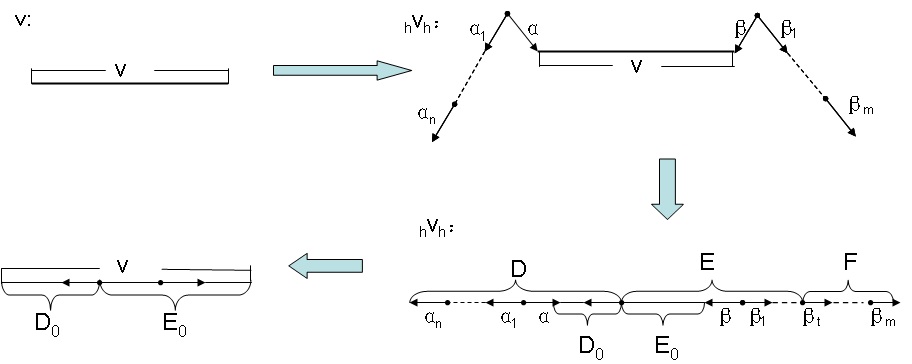}
\end{center}
    If $E^{-1}=E'$, there is a two-sided ad-pair
    $$(D_0, E_0, 1_{e(\bz)})\times (D'\bz^{-1}_{t}\cdots\bz^{-1}_1\bz, E^{-1}_0,F')$$
    $$\in\ff(v)\times \fs(w).$$
\end{itemize}
 \item[II.] $l(D)< l(V_\az\az).$ Let $V_\az=\az^{-1}_n\cdots\az^{-1}_2\az^{-1}_1$, assume that
  $D=\az^{-1}_n\cdots\az^{-1}_{j+1}$ where $0\leq j\leq n$ with $\az^{-1}_{n+1}:=1_{e(\az_n)}$ and $E=\az^{-1}_{j}\cdots\az^{-1}_1\az E_0$ with $\az^{-1}_0:=1_{s(\az)}$. Then $l(E_0)\geq 0$ since $l(DE)>l(V_\az)$. Similarly as in I, the definition of factor string implies that $l(F)\leq l(\az{v_h}).$
  If $l(F)=l(\az{v_h})$, then there is nothing to prove. We assume $l(F)<l(\az{v_h})$ and consider the following two subcases:
  \begin{itemize}
    \item[(II.1)]$l(\az{v_h})>l(F)>l(\bz^{-1}V^{-1}_\bz).$ We rewrite $F=F_0\bz^{-1}V^{-1}_\bz$ with $l(F_0)>0.$ Then there exists a two-sided ad-pair
   $$(1_{e(\az)}, E_0, F_0)\times (D'\az^{-1}_{j}\cdots\az^{-1}_1\az,E_0,F')$$
   $$\in \ff(v)\times \fs(w)$$
   if $E'=E$, since $l(D'\az^{-1}_{j}\cdots\az^{-1}_1\az)>0$ and $l(F_0)>0.$  Similarly, if $E'=E^{-1}$,
   the following ad-pair
   $$(1_{e(\az)}, E_0, F_0)\times (D',E^{-1}_0,\az^{-1}\az_1\cdots\az_{j}F')$$
   is two-sided in $\ff(v)\times \fs(w).$ See the following figure:
\begin{center}
\includegraphics[height=1.7in]{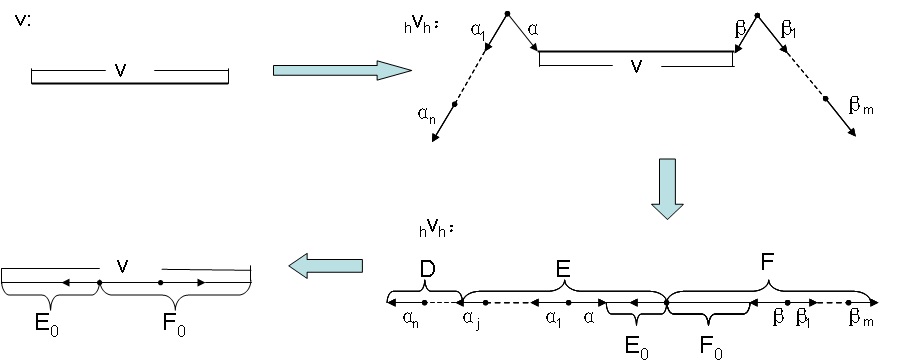}
\end{center}
    \item[(II.2)]$l(F)\leq l(V^{-1}_\bz).$ Let $V^{-1}_\bz=\bz_1\cdots\bz_m$ with $m\geq 0$ where $\bz_0:=1_{s(\bz)}$, and assume $F=\bz_{t}\cdots\bz_m$ with $0\leq t\leq m+1$ where $\bz_{m+1}:=e(\bz_m).$ Therefore, we can rewrite
        $$E=\az^{-1}_{j}\cdots\az^{-1}_1\az v\bz^{-1}\bz_1\cdots\bz_{t-1}.$$
        Then we obtain a two-sided ad-pair
        $$(1_{e(\az)},v,1_{e(\bz)})\times (D'\az^{-1}_{j}\cdots\az^{-1}_1\az, v,\bz^{-1}\bz_1\cdots\bz_{t-1}F')$$
                     $$=(1_{e(\az)},v,1_{e(\bz)})\times ({D"},v,F")\in\ff(v)\times \fs(w)$$
    when $E=E'$. See the following figure
\begin{center}
\includegraphics[height=2.5in]{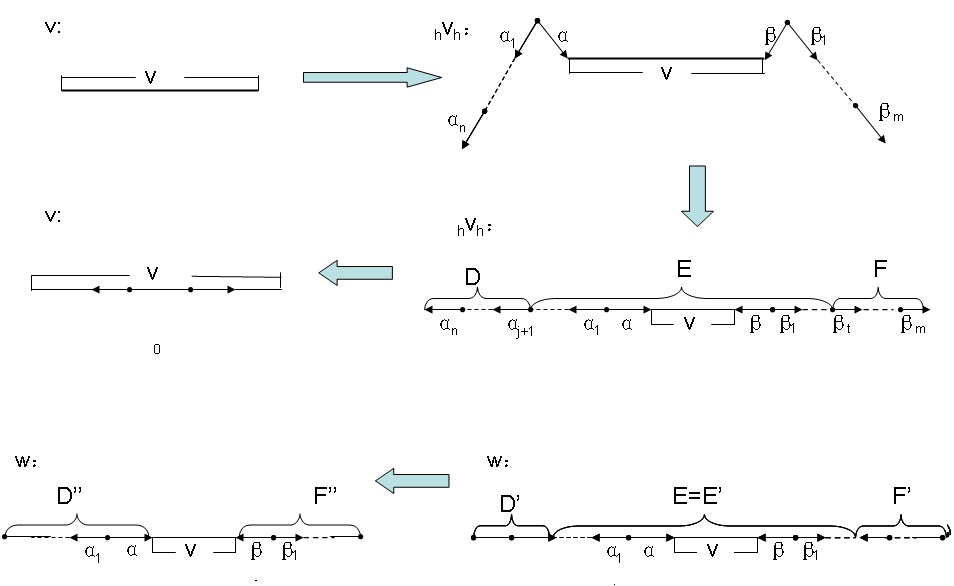}\label{fig}
\end{center}
 If $E^{-1}=E',$ we get the following two-sided ad-pair
$$(1_{e(\bz)},v,1_{e(\bz)})\times (D'\bz^{-1}_{t-1}\cdots\bz^{-1}_1\bz, v^{-1},\az^{-1}\az_1\cdots\az_{j}F')$$
                     $$\in\ff(v)\times \fs(w).$$
  \end{itemize}
Therefore we complete the proof of the claim.
By using Similar discussion as in Lemma \ref{case1}, Lemma \ref{lem-important} implies that $f_p$ factors through projective module. It completes the proof.
  \end{itemize}
\end{pf}

\begin{thm}Let $w,v$ be two strings over $A$. Then $\Ext^1_A(M(w),M(v))\neq 0$ if and only if $w,v$ admit extensions
from $w$ to $v.$
\end{thm}
\begin{pf} By Corollary \ref{cor-imp}, it suffices to show that if $w,v$ do not have extensions from $w$ to $v,$ then $\Ext^1_A(M(w),M(v))=0.$

By Lemma \ref{case1}, Corollary \ref{case2}  and Lemma \ref{case3}, we only need to consider the case that $v$  both starts and ends on a peak.
It suffices to show that $\Hom_A(\t^{-1}M(v), M(w))=0$ by Auslander-Reiten formula.

Assume $v$ both starts and ends on a peak, Lemma \ref{lem-inj} implies that we only need to consider non-injective module $M(v)$. Therefore, there exists $\az, \bz$ in $Q_1$ such that we can write
$$v=\az_n\cdots\az_1\az^{-1}v_0\bz\bz^{-1}_1\cdots\bz^{-1}_m=U^{-1}_\az\az^{-1}v_0\bz U_\bz$$ with $m,n\geq 0$.
By Theorem \ref{BR87}, we have
$$\t^{-1}M(v)=M({_cv}_c)=M(v_0).$$

Suppose $\Hom_A(\t^{-1}M(v), M(w))=\Hom_A(M(v_0), M(w))\neq 0$ which means by Theorem \ref{C-B89} that $v_0,w$ have an ad-pair
$(D,E,F)\times (D',E',F')\in\ff(v_0)\times \fs(w).$  It implies that
 $$(U^{-1}_\az\az^{-1}D,E,F\bz U_\bz)\in \ff(v).$$
\begin{center}
\includegraphics[height=2in]{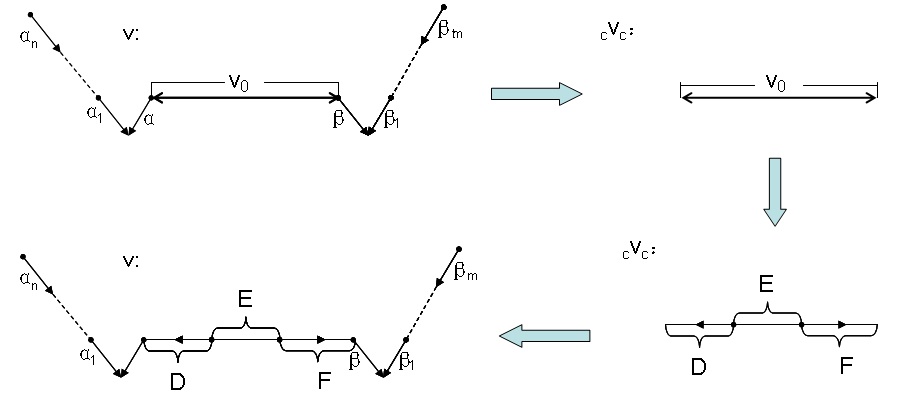}\label{fig2}
\end{center}
Therefore there is a two-sided ad-pair
 $$(U^{-1}_\az\az^{-1}D,E,F\bz U_\bz)\times (D',E',F')\in \ff(v)\times \fs(w)$$
which contradicts the hypothesis that $w,v$ do not have extensions from $w$ to $v$. Then
$\Hom_A(\t^{-1}M(v), M(w))=0$ which implies $\Ext^1_A(M(w),M(v))=0.$
\end{pf}

\section*{Acknowledgments}
The author is indebted to professor Feng Luo for pointing out Theorem \ref{mocher} when he was giving series of lectures
in Mathematical center in Tsinghua University. He is very grateful for professor Thomas Br\"ustle, professor Bin Zhu and Yu Zhou
with many helpful discussions and he also thanks professor Jan Schr\"oer for many helpful comments and for telling him that a conjecture in a previous
version of the paper is wrong.


\begin{thebibliography}{99}
\bibitem[A09]{A09} Amiot C., Cluster categories for algebras of global dimension 2 and quivers with potential. {\em Annales de l'Institut Fourier} {\bf 59} (2009): 2525-2590.

\bibitem[ABCP10]{ABCP10} Assem I., Br\"ustle T., Charbonneau-Jodoin G. and P-G. Plamondon. Gentle algebras arising from surface triangulations,
{\em J. Algebra and Number Theory} {\bf 4} (2010): 201-229.

\bibitem[AD12]{AD12}Assem I. and Dupont G., Moduls over cluster-tilted algebras determined by their dimension vectors, http://arxiv.org/abs/1202.5698.

\bibitem[AS87]{AS87}Assem I. and  Skowro\'nski, A., Iterated tilted algebras of type $\widetilde{A_n}$, {\em Math. Z}. {\bf 195} (1987): 269-290.

\bibitem[BMRRT06]{BMRRT06} Buan A. B.,Marsh R., Reineke M., Reiten I. and Todorov G., Tilting theory
and cluster combinatorics, {\em Adv. Math}. {\bf 204} (2006), 572-612.

\bibitem[BR87]{BR87} Butler M. C. R. and Ringel C. M., Auslander-Reiten sequences with few
middle terms. {\em Comm. in Algebra} {\textbf{15}} (1987):145-179.

\bibitem[BZ11]{BZ11}Br\"ustle T. and Zhang J., On the cluster category of a marked surface without punctures. {\em J. Algebra and Number Theory}  {\bf 5-4} (2011):529--566.

\bibitem[C11]{C11}Carroll A. T., Generic Modules for string algebras http://arxiv.org/abs/1111.5064.

\bibitem[C-B89]{C-B89}
W. W. Crawley-Boevey., Maps between representations of zero-relation algebras
{\em J. Algebra} \textbf{126}(2), (1989): 259-263.


\bibitem[DWZ08]{DWZ08} Derksen H., Weyman J. and Zelevinsky A., Quivers with potentials and
their representations I: Mutations. {\em Selecta Math., New Series} {\bf 14} (2008):59-119.

\bibitem[FST08]{FST08} Fomin S., Shapiro M. and Thurston D.,Cluster algebras and triangulated
surfaces. Part I: Cluster complexes. {\em Acta Mathematica} {\bf 201} (2008):83-146.

\bibitem[FZ02]{FZ02} Fomin S. and Zelevinsky A.,
 Cluster algebras I. Foundations. \emph{ J. Amer. Math. Soc. }
{\bf 15}(2) (2002):497--529 (electronic).

\bibitem[G72]{G72} Gabriel P., Unzerlegbare Darstellungen. I. {\em Manuscripta Math.}, 6: 71-103; correction, ibid. {\bf 6} (1972), 309,
1972.

\bibitem[GP12]{GP12} Geng S. F. and Peng L. G., The dimension vectors of indecomposable modules of cluster-tilted algebras and the Fomin-Zelevinsky denominators conjecture. \emph{Acta Mathematica Sinica, English Series,} 28(3), (2012): 581-586.


\bibitem[K10]{K10} Keller B., Cluster algebras, quiver representations and triangulated categories. In Triangulated
categories, {\em LMS Lecture Notes Ser.}, {\bf 375}, Pages 76-160. Cambridge Univ. Press, Cambrdge, 2010.

\bibitem[M83]{M83} Mosher L., Pseudo-Anosovs on punctured surfaces, Ph.D. thesis, Princeton
University, 1983.

\bibitem[LF09]{LF09} Labardini-Fragoso D.,  Quivers with potentials associated to triangulated
surfaces. {\em Proc. London Math. Soc.}  {\bf 98} (2009):797-839.

\bibitem[S99]{S99} J. Schr\"oer., Modules without self-extensions over gentle algebras.
{\em J. Algebra} \textbf{216}, (1999):178-189.

\bibitem[T10]{T10}THURSTON D., Geometric intersection of curves on surfaces. Preprint
http://www.math.columbia.edu/~ dpt/DehnCoordinates.ps

\bibitem[WW85]{WW85} Wald B. and Waschb\"usch J., Tame biserial algebras,{\em J. Algebra} {\bf 95} no. 2, (1985):480-500.

\end{thebibliography}
\end{document}